\documentclass{amsart}
\usepackage{amssymb, amsmath, amsfonts, amscd}

\input xypic

\theoremstyle{plain}
\newtheorem{theorem}{Theorem}[section]

\newtheorem{lemma}[theorem]{Lemma}
\newtheorem{proposition}[theorem]{Proposition}

\newtheorem{example}[theorem]{Example}

\theoremstyle{definition}
\newtheorem{definition}[theorem]{Definition}

\theoremstyle{remark}

\numberwithin{equation}{theorem}

\newcommand{\m}[1]{\mathcal{#1}}

\renewcommand{\P}{\mathbf{P} }

\renewcommand{\O}{\mathcal{O}}

\renewcommand{\L}{\mathcal{L}}

\newcommand{\End}{\operatorname{End} }

\renewcommand{\H}{\operatorname{H} }
\newcommand{\SL}{\operatorname{SL} }
\newcommand{\U}{\operatorname{U}}

\newcommand{\Sym}{\operatorname{Sym} }

\renewcommand{\lg}{\mathfrak{g}}
\renewcommand{\sl}{\mathfrak{sl} }

\renewcommand{\ln}{\mathfrak{n}}
\newcommand{\lb}{\mathfrak{b}}
\newcommand{\lp}{\mathfrak{p}}
\newcommand{\lpv}{\mathfrak{p}(V_\bullet)}
\newcommand{\lnv}{\mathfrak{n}(V_\bullet)}
\newcommand{\lnvop}{\mathfrak{n}(V_\bullet)_{op}}
\newcommand{\lpvop}{\mathfrak{p}(V_\bullet)_{op}}
\newcommand{\lqv}{\mathfrak{q}(V_\bullet)}

\newcommand{\llqv}{\tilde{\mathfrak{q}}(V_\bullet)}

\newcommand{\llnv}{\tilde{\mathfrak{n}}(V_\bullet)}

\newcommand{\lh}{\mathfrak{h}}
\newcommand{\lmv}{\mathfrak{m}(V_\bullet)}

\newcommand{\oneg}{\mathbf{1}}

\renewcommand{\m}{\mathfrak{m} }

\begin{document}

\title{On the canonical filtration of an irreducible representation}

\author{Helge Maakestad }
\address{Hoegskolen i Bergen} 
\email{h\_maakestad@hotmail.com}



\date{May 2012} 

\begin{abstract} The aim of this paper is to study the canonical filtration 
$L_l(\lambda)$ 
of an irreducible finite dimensional $\sl(V)$-module
$L(\lambda)$ using the universal enveloping algebra $\U(\sl(V))$ and the annihilator ideal $ann(v)$ of a highest weight vector $v$
in $L(\lambda)$. We give a basis for 
$L_l(\lambda)$ 
in terms of the universal enveloping algebra of the nilpotent radical 
of an opposite parabolic sub algebra to the stabilizer Lie algebra of a flag $V_\bullet$ in $V$. We also calculate the dimension of 
$L_l(\lambda)$ 
as a function of $l$.
\end{abstract}

\maketitle
\tableofcontents

\section{Introduction}

The aim of this paper is to study the canonical filtration $L_l(\lambda)$ of an irreduible finite dimensional $\SL(V)$-module $L(\lambda)$ 
using the universal enveloping algebra $\U(\sl(V))$ and the annihilator ideal $ann(v)$ of a highets weight vector $v$ in $L(\lambda)$.
The stabilizer Lie algebra $\lp_v$ of the vector $v$ may be realized as the stabilizer Lie algebra $\lpv$ of a flag $V_\bullet$ in $V$.
We may construct the opposite parabolic subalgebra $\lpvop$ of $\lpv$ relative to a choice of roots $R$ of $\sl(V)$. The nilpotent radical
$\lnvop$ of $\lpvop$ is a sub Lie algebra of $\sl(V)$ such that $\sl(V)=\lnvop \oplus \lpv$. In the main theorem of the paper 
(see Theorem \ref{maintheorem}) we use the universal enveloping algebra
$\U(\lnvop)$ and its canonical filtration $\U_l(\lnvop)$ to construct a basis for the $\lpv$-module $L_l(\lambda)$. We also calculate the dimension
of $L_l(\lambda)$ as a function of $l$.
This result is applied in the paper  \cite{maakestad0} to give a geometric construction of the terms in the canonical filtration
in terms of jet bundles of line bundles on the flag variety $\SL(V)/P$. The aim of the construction is to apply this to the study of discriminants of
linear systems on flag varieties (see Example \ref{discriminants} and the preprint \cite{maakestad2}).

\section{Root systems, flags and nilpotent radicals of parabolic Lie algebras}

Let $K$ be an algebraically closed field of characteristic zero and let $V=K^N$ be $N$-dimensional vector space over $K$.
Let $l_1,\ldots , \l_{k+1}\geq 1$ be integers with $l_1+\cdots +l_{k+1}=N$ and let $d_i=l_1+\cdots +l_i$ for $i=1,\ldots, k+1$.
Let 
\[ V_\bullet: 0\neq V_1\subsetneq V_2\subsetneq \cdots \subsetneq V_k\subsetneq V_{k+1}=V \]
be a flag of vector spaces in $V$ with $dim_K(V_i)=d_i$. Let $\lg=\sl(V)$ be the special Lie algebra on $V$ and let
$\lp(V_\bullet) \subseteq \lg$ be the Lie algebra of elements fixing the flag $V_\bullet$. By this we mean the following:
$\lp(V_\bullet)$ is the set of elements $x\in \lg$ such that for any $i=1,\ldots , k$ and any vector $v\in V_i$ it follows 
$x(v)\in V_i$. One checks $\lp(V_\bullet)$ is a sub Lie-algebra of $\lg$.

\begin{definition} Let $\lp(V_\bullet)$ be the \emph{stabilizer Lie algebra} of the flag $V_\bullet$.
\end{definition}

Let $\lh\subseteq \lg$ be the abelian Lie algebra of
diagonal matrices with trace equal to zero. It follows from \cite{dixmier} the pair $(\lg, \lh)$ is a split semi-simple Lie algebra.
There is the adjoint representation $ad:\lh \rightarrow \End(\lg)$ and this determines the \emph{roots of $\lg$ relative to $\lh$} denoted
$R=R(\lg,\lh)$. It follows $R$ is the set of elements $\lambda\in \lh^*$ with $\lg_\lambda \neq 0$
(see \cite{dixmier}, Section 1.2.13 for the definition of $\lg_\lambda)$).
Let $E_{ij}$ be the $N\times N$-matrix with a $1$ on place $(i,j)$ and zero's elsewhere. 
Let
\[ \lh^*=K\{L_1,\ldots ,L_N\}/L_1+\cdots +L_N \]
with 
\[ L_i(E_{jj})=\delta_{ij}.\]
It follows
\[ R=\{L_i-L_j : 1\leq i \leq j \leq N\}.\]
Let $\lambda_{ij}=L_i-L_j$.
It follows $\lg_{\lambda_{ij}}=KE_{ij}$ is the one dimensional vector space over $K$ with $E_{ij}$ as a basis.
It follows    $\lpv$    consists of matrices on the form
\[
x=
\begin{pmatrix} A_1 & * & \cdots & *  \\
                0 & A_2 & \cdots & * \\
              \vdots &  \cdots &  \cdots & \vdots \\
                0 &     \cdots &  \cdots & A_{k+1}
\end{pmatrix}
\]
where $tr(x)=0$ and $A_i$ is an $l_i\times l_i$-matrix with coefficients in $K$. 
Let
\[ I_{l}=\{ (i,j) : d_{l-1}<i\leq d_l\text{ and }d_{l-1}<j \leq N\}\]
for $l=1,\ldots, k+1$. Let $I=\cup_{l=1}^{k+1}I_{l}$.

\begin{lemma} It follows $\lpv= \oplus_{(i,j)\in I}\lg_{\lambda_{ij}}$.
\end{lemma}
\begin{proof} The proof is left to the reader as an exercise.
\end{proof}

\begin{definition} Let $\lpvop=\oplus_{(i,j)\in I}\lg_{\lambda_{ji}}$ be the \emph{opposite Lie algebra of} $\lpv$.
Let $\lnvop$ be the nilpotent radical of $\lpvop$.
\end{definition} 

Note: The Lie algebras $\lpvop$ and $\lnvop$ depend on a choice of maximal abelian sub algebra $\lh$ in $\lg$.

Let $\lnv \subseteq \lpv$ be the subalgebra  consisting of elements with $A_i$ the zero matrix for $i=1,\ldots k+1$. 
Let $\lqv$ be the subalgebra with zero's off the diagonal matrices $A_i$ for $i=1,\ldots k+1$. It follows there is a direct sum decomposition
$\lpv=\lqv\oplus \lnv$. 

\begin{lemma} The following holds: $\lpv \cong \lnv \times \lqv$ is a semi direct product of Lie algebras.
\end{lemma}
\begin{proof} The proof is left to the reader as an exercise.
\end{proof}

By definition it follows $dim_K(\lqv)=l_1^2+\cdots +l_{k+1}^2-1$.
Let $\llqv$ be the subalgebra of $\lqv$ consisting of matrices with $tr(A_i)=0$ for $i=1,\ldots ,k+1$. It follows
$\llqv=\sl(l_1)\oplus \cdots \oplus \sl(l_{k+1})$. By definition $dim_K(\llqv)=l_1^2+\cdots +l_{k+1}^2-k-1$. It follows the codimension 
of $\llqv$ in $\lqv$ is $k$. Let $D_i$ be the matrix with $A_i=l_{i+1}I_{l_i}$ and $A_{i+1}=-l_iI_{i+1}$ on 
the diagonal. $A_j=0$ for $j\neq i,i+1$, zero's elsewhere.
Here $I_j$ is the $j\times j$ identity matrix. It follows $tr(D_i)=0$ for $i=1,\ldots k$.

\begin{lemma} \label{dimension}  The set $D=\{D_1,\ldots ,D_k\}$ is linearly independent over $K$.
\end{lemma}
\begin{proof} The proof is left to the reader as an exercise.
\end{proof}

Let $L_D$ be the subalgebra of $\lqv$ spanned by $D$ over $K$. It follows $L_D$ is an abelian Lie algebra of dimension $k$. One checks
$[x, d]=0$ for any $x\in \lpv$ and $d\in L_D$. There is a map
\[ \phi: \llqv\oplus L_D \rightarrow \lqv \]
defined by
\[ \phi(x,\omega)=x+\omega .\]

\begin{lemma} The map $\phi$ is an isomorphism of Lie algebras.
\end{lemma}
\begin{proof} Since $\llqv\cap L_D=\{0\}$ it follows $\phi$ is an injection. By Lemma \ref{dimension} it follows $\phi$ is an isomorphism
of vector spaces. We get the following calculation:
\[ \phi([(x,u),(y,v)])=\phi([x,y],0)=[x,y].\]
We get
\[ [\phi(x,u),\phi(y,v)]=[x+u,y+v]=[x,y]+[x,v]+[u,y]+[u,v]=[x,y].\]
It follows $\phi$ is an isomorphism of Lie algebras.
\end{proof}

It follows 
\[ \lqv\cong \sl(l_1)\oplus \cdots \oplus \sl(l_{k+1})\oplus L_D \]
is the direct sum of a semi simple Lie algebra and an abelian Lie algebra hence it is linearly reductive.
We get a direct sum decomposition $\lg=\llqv \oplus L_D \oplus \lnv$ of vector spaces.
Let $\llnv=L_D\oplus \lnv$. Since $\lnv$ is a solvable Lie algebra and $L_D$ is abelian hence solvable it follows $\llnv$ is a solvable
Lie algebra. 

Let $D$ be the following matrix:
\[
D=
\begin{pmatrix} D_1     &          0          &  \cdots      &     0      \\
                 0      &     D_2             &  \cdots      &      0     \\
                 \vdots &    \vdots           &  \cdots      &     \vdots \\
                 0      &    0                &  \cdots      &   D_{k+1}  
\end{pmatrix},
\]
where $D_i=\alpha_iI_{l_i}$ where $\alpha_i\in K$ and $I_{l_i}$ is the $l_i\times l_i$ identity matrix.
Assume $tr(D)=\alpha_1l_1+\cdots +\alpha_{k+1}l_{k+1}=0$. It follows
\[ \alpha_{k+1}=\frac{-\alpha_1 l_1-\cdots -\alpha_k l_k}{l_{k+1}}.\]
We get an element $D\in \lg$.
Let $1 < i < j < N$ and  assume $x$ is the following matrix:
\[
x=
\begin{pmatrix}   E_1          &   \cdots  &      0       &    \cdots   &   0        \\
                   0           &   \cdots  &      0       &    \cdots   &   0        \\
                  \vdots       &   \cdots  &      U_{ij}  &    \cdots   &   0        \\
                  \vdots       &   \cdots  &     0        &    \cdots   &   0        \\
                  0            &  \cdots   &      0       &   \cdots    &  E_{k+1}  
\end{pmatrix},
\]
where $tr(x)=0$ and in the $(i,j)$-block there is an $l_i\times l_j$-matrix $U_{ij}$. It follows $x\in \lg$.

\begin{lemma} \label{commutatorone} Assume $\alpha_j=\alpha_i-1$. It follows $[D,x]$ is the matrix with $U_{ij}$ of the $(i,j)$-block and 
zeros elsewhere.
\end{lemma}
\begin{proof} One checks that the matrix $[D,x]$ is the matrix with $D_iU_{ij}-U_{ij}D_j$ on the $(i,j)$-bloc. 
Since $D_i=\alpha_iI_{l_i}$ it follows $D_iU_{ij}=\alpha_iU_{ij}$. We get $D_iU_{ij}-U_{ij}D_j=(\alpha_i-\alpha_j)U_{ij}=U_{ij}$.
The lemma is proved.
\end{proof}

\begin{lemma} \label{commutator} There is an equality $[\lpv, \lpv]=\sl(l_1)\oplus \cdots \oplus \sl(l_{k+1})\oplus \lnv$ of Lie algebras.
\end{lemma}
\begin{proof} Let $x,y \in \lpv$. It follows $[x,y]$ has commutators $[A_i,B_i]$ where $A_i,B_i$ are $l_i\times l_i$-matrices
on the block diagonals. It follows $[x,y]\in \sl(l_1)\oplus \cdots \oplus \sl(l_{k+1})\oplus \lnv$. By Lemma \ref{commutatorone}
it follows $\lnv \subseteq [\lpv,\lpv]$. Since $[\sl(n),\sl(n)]=\sl(n)$ it follows
\[ \sl(l_1)\oplus \cdots \oplus \sl(l_{k+1})\subseteq [\lpv, \lpv].\]
The claim of the lemma is proved.
\end{proof}

\begin{proposition} \label{nilpotent} The ideal $\llnv$ is the radical of $\lpv$. The ideal $\lnv$ is the nilpotent radical of $\lpv$.
\end{proposition}
\begin{proof} Assume there is a strict inclusion $\llnv \subseteq I \subseteq \lpv$ of solvable ideals in $\lpv$. 
Via the quotient morphism $\lpv \rightarrow \lpv/\llnv$ we get a non-zero solvable ideal 
$J\subseteq \lpv/\llnv\cong \sl(l_1)\oplus \cdots \oplus \sl(l_{k+1})$. Since $\lpv/\llnv$ is a semi simple Lie algebra it contains no 
non-trivial solvable ideals. This is a contradiction hence $\llnv$ is the maximal solvable ideal in $\lpv$ and the first part of the 
 proposition is proved. By definition it follows the nilpotent radical of $\lpv$ equals $[\lpv,\lpv]\cap \llnv$. By Lemma \ref{commutator}
there is an equality of Lie algebras
\[ [\lpv, \lpv]=\sl(l_1)\oplus \cdots \oplus \sl(l_{k+1})\oplus \lnv.\]
It follows
\[ [\lpv,\lpv]\cap \llnv =\sl(l_1)\oplus \cdots \sl(l_{k+1})\oplus \lnv\cap L_D\oplus \lnv=\lnv.\]
The proposition is proved.
\end{proof}

Let $\lmv$ be the subalgebra of $\lg$ consisting of matrices on the form
\[
x=
\begin{pmatrix} A_1      &  0        & \cdots   &   \cdots    &  0     \\ 
                *        &  A_2      & \cdots   & \cdots      &  0     \\ 
                \vdots   & \cdots    & \cdots   & \cdots      &  0     \\
               *         & \cdots    & \cdots   &   *         & A_{k+1}    
\end{pmatrix}
\]
where $A_i$ is the $l_i\times l_i$-matrix with zero's as entries. It follows $\lg=\lmv\oplus \lpv$ as vector space.

\begin{proposition} \label{dimension} There is an equality $\ln(V_\bullet)_{op}=\lmv$ of Lie algebras.
The following formula holds:
\[ dim_K(\ln(V_\bullet)_{op})=\sum_{1\leq j<i \leq k+1}l_il_j.\]
\end{proposition}
\begin{proof} The proof of the lemma follows from Proposition \ref{nilpotent}.
\end{proof}

It follows there is a direct sum decompositon
\[ \sl(E))=\ln(V_\bullet)_{op}\oplus \lp(V_\bullet)   \]
of vector spaces for every flag $V_\bullet$ in $E$. 

\begin{definition}
Let the Lie algebra $\ln(V_\bullet)_{op}$ be the \emph{opposite radical} of the flag
$V_\bullet$.
\end{definition}

\section{On the canonical filtration of an irreducible representation}

Let $L(\lambda)$ be an irreducible finite dimensional  left $\lg$-module with highest weight vector $v$ and highest weight $\lambda$. 
The line $L_v$ defined by $v$ is a unique line in $L(\lambda)$. Let $\lp_v$ be the subalgebra
of $\lg$ consisting of elements fixing the line $L_v$. It follows $\lp_v$ is a parabolic subalgebra of $\lg$. Moreover it follows
from Theorem \ref{irreducible} that $\lp_v=\lp(V_\bullet)$ for some flag $V_\bullet$ in $V$.
By the previous section it follows there is a decomposition $\lg=\ln(V_\bullet)_{op}\oplus \lp(V_\bullet)$
where $\lnvop$
is the nilpotent radical of an opposite parabolic algebra $\lp(V_\bullet)_{op}$ of $\lp(V_\bullet)$.
We may choose a basis $y_1,\ldots ,y_{k-1}, y_k$ for $\lpv$ satisfying the following property:
There is an exact sequence
\[ 0\rightarrow ker(\rho) \rightarrow \lpv \rightarrow^{\rho}  \End(L_v)\rightarrow 0\]
and the set $y_1,\ldots ,y_{k-1}$ is a basis for the Lie-algebra $Ker(\rho)$. Moreover $y_k(v)=v$. We may choose a basis $z_1,\ldots ,z_l$ 
for $\lnvop$.
In the following we consider the universal enveloping algebra $\U(\lg)$ and its unit element $\oneg$.
Let 
\[ B_j=\{ z_1^{p_1}\cdots z_l^{p_l}y_j^{q_j+1}\cdots y_k^{q_k} : p_i,q_j\geq 0\} \]
for $j=1,\ldots ,k-1$.
Let 
\[ B_k=\{z_1^{p_1}\cdots z_l^{p_l}y_k^{q_k}(y_k-\oneg) : p_i,q_k\geq 0\} .\]

\begin{lemma} \label{basis1} The set $B_i$ is a linearly independent set of elements in $\U(\lg)$ over the field $K$ for $i=1,\ldots , k$.
\end{lemma}
\begin{proof} Consider the set $B_i$ for $i=1,\ldots , k-1$. This is a lineary independent set of vectors by the PBW-theorem.
Consider the set $B_k$. Assume there is a relation
\[ \sum_{p,q_k} \alpha_{p,q_k}z_1^{p_1}\cdots z_l^{p_l}y_k^{q_k}(y_k-\oneg)=0 \]
It follows 
\[ (\sum_{p,q_k} \alpha_{p,q_k}z_1^{p_1}\cdots z_l^{p_l}y_k^{q_k})(y_k-\oneg)=\omega (y_k-\oneg)=0.\]
Since $\U(\lg)$ is an integral domain it follows $\omega=0$. The lemma is proved.
\end{proof}

\begin{lemma} \label{basis2} Let $1\leq j < i \leq k$. The following holds: $B_j\cap B_i=\emptyset$.
\end{lemma}
\begin{proof} Assume $1< i \leq k-1$ and let $\omega=z_1^{p_1}\cdots z_l^{p_l}y_i^{p_i+1}\cdots y_k^{p_k}\in B_i$. If
$\omega \in B_j$ it follows 
\[ \omega=z_1^{p_1}\cdots z_l^{p_l}y_j^{p_j+1}\cdots y_k^{p_k} .\]
Since $\omega \in B_i$ it follows $p_j+1=0$ hence $p_j=-1$ which is a contradiction. Hence $B_j\cap B_i=\emptyset$.
Assume $i=k$ and $1\leq j < i$. Let $\omega=z_1^{p_1}\cdots z_l^{p_l}y_k^{p_k}(y_k-\oneg)\in B_k$. If $\omega \in B_i$ it follows
\[ \omega= z_1^{p_1}\cdots z_l^{p_l}y_i^{p_i+1}\cdots y_k^{p_k} \]
with $i<k$. Hence $p_i+1=0$ and $p_i=-1$ which is a contradiction. The lemma is proved.
\end{proof}

Let $B=\{B_1,\ldots , B_k\}$. It follows from Lemma \ref{basis1} and \ref{basis2} that $B$ is a linearly independent set of 
vectors in $\U(\lg)$.

Let $char(\rho) = \U(\lg)\{ x-\rho(x)\oneg : x\in \lpv\}$ be the \emph{character ideal} of the character $\rho$ of $\lpv$.

By definition  $char(\rho)$ is a left ideal in $\U(\lg)$. The aim of this section is to prove that $B$ is a basis for $char(\rho)$ 
over the field $K$. Let $char_l(\rho)=char(\rho)\cap \U_l(\lg)$ for any integer $l\geq 1$. Here $\U_l(\lg)\subseteq \U(\lg)$ is the canonical
filtration of $\U(\lg)$. Let $K\{B\}$ be the subspace of $\U(\lg)$ spanned by the set $B$. Let $K\{B\}_{r+1}=K\{B\}\cap \U_{r+1}(\lg)$.

\begin{lemma} \label{basis3} The following formula holds for any integer $k\geq 1$:
\[ dim_K (K\{B\}_{r+1} )=\sum_{i=0}^{k-1} \binom{r+l+k-i}{l+k-i}.\]
\end{lemma}
\begin{proof} Let $1\leq i \leq k-1$ and let 
\[ B_i^{r+1}=\{ z_1^{p_1}\cdots z_l^{p_l}y_j^{q_j+1}\cdots y_k^{q_k}:0\leq \sum p_j+\sum q_j\leq r\}.\]
Let
\[ B_k^{r+1}=\{z_1^{p_1}\cdots z_l^{p_l}y_k^{q_k}(y_k-\oneg): 0\leq \sum p_j+q_k \leq r\} .\]
It follows $B_i^{r+1}, B_k^{r+1}\subseteq \U_{r+1}(\lg)$. One calculates that the number of elements in 
$B_i^{r+1}$ equals $\binom{r+l+k-i+1}{l+k-i+1}$. Also the number of elements in $B_k^{r+1}$ equals $\binom{l+r+1}{l+1}$.
We get the formula
\[ dim_K( K\{B\}_{r+1})=\sum_{i=1}^{k-1}\binom{r+l+k-i+1}{l+k-i+1}+\binom{l+r+1}{l+1}.\]
We change index and let $j=i-1$ to get
\[ dim_K( K\{B\}_{r+1})=\sum_{j=0}^{k-2}\binom{r+l+k-j}{l+k-j}+\binom{l+r+1}{l+1}=\]
\[ \sum_{j=0}^{k-1}\binom{r+l+k-i}{l+k-i}.\]
The lemma is proved.
\end{proof}

\begin{lemma} \label{basis4} Let $k\geq 1$ be an integer. The following formula holds:
\[ \sum_{i=0}^{k-1}\binom{r+l+k-i}{l+k-i}+\binom{r+l+1}{l} =\binom{r+l+k+1}{l+k}  .\]
\end{lemma}
\begin{proof} The formula is proved by induction.
\end{proof}

By definition $dim_K(\lg)=N^2-1$, $dim_K(\lp)=k$ and $dim_K(\ln)=l$. It follows $N^2-1=l+k$ and
$N^2=l+k+1$. We get
\[ dim_K(\U_{r+1}(\lg)=\binom{N^2+r}{N^2-1}=\binom{r+l+k+1}{l+k} \]
and
\[ dim_K(\U_{r+1}(\ln)=\binom{r+l+1}{l}.\]
By definition it follows there is an inclusion $K\{B\}\subseteq char(\rho)$ and $K\{B\}_{r+1}\subseteq char_{r+1}(\rho)$ of
vector spaces. 

\begin{proposition} \label{directsum} There is for every integer $r\geq 0$ a direct sum decomposition
\[ \U_{r+1}(\lg)=\U_{r+1}(\lnvop)\oplus char_{r+1}(\rho)\]
of vector spaces.
\end{proposition}
\begin{proof} By Lemma \ref{basis3} the set $B$ spans the sub vector space $K\{B\}_{r+1}$ of $char_{r+1}(\rho)$ of dimension
\[ \sum_{i=0}^{k-1}\binom{r+l+k-i}{l+k-i}.\]
Since $\U_{r+1}(\lnvop)\cap K\{B\}_{r+1}=0$ it follows the sum
\[ \U_{r+1}(\lnvop)+K\{B\}_{r+1} \]
is direct.
There is an inclusion
\[ \U_{r+1}(\lnvop)\oplus K\{B\}_{r+1}\subseteq \U_{r+1}(\lnvop)+char_{r+1}(\rho) \subseteq \U_{r+1}(\lg) .\]
By Lemma \ref{basis3} and \ref{basis4} it follows $dim_K(\U_{r+1}(\lnvop)\oplus K\{B\}_{r+1})=dim_K(\U_{r+1}(\lg))$ It follows
\[ K\{B\}_{r+1}=char_{r+1}(\rho).\]
The proposition is proved.
\end{proof}

\begin{definition} Let $M(\lambda)=\U(\lg)\otimes_{\U(\lp)}L_v$ be the \emph{generalized Verma module} associated to $\lambda$.
Let $M_l(\lambda)=\U_l(\lg)\otimes_{\U(\lp)}L_v$ be its \emph{canonical filtration}.
\end{definition}

A generalized Verma module $M(\lambda)$ associated to a character $L_v$ of a parabolic Lie algebra $\lp$ is sometimes called a
\emph{scalar generalized Verma module}.

There is a unique maximal left $\lg$-sub module $K(\lambda)\subseteq M(\lambda)$ and an exact sequence of $\lg$-modules
\[ 0\rightarrow K(\lambda) \rightarrow M(\lambda) \rightarrow L(\lambda) \rightarrow 0 .\]
By definition it follows
\[ K(\lambda)=\{x\otimes v : x(v)=0 \} ,\]
hence it follows
\[ K(\lambda)=ann(v)\otimes L_v \]
where $ann(v)\subseteq \U(\lg)$ is the annihilator ideal of the highest weight vector $v$.

We get for every integer $l\geq 1$ exact sequences
\[ 0\rightarrow ann_l(v)\otimes L_v \rightarrow M_l(\lambda) \rightarrow L_l(\lambda) \rightarrow 0 .\]

\begin{definition} Let $\{L_l(\lambda)\}_{l\geq 1}$ be the \emph{canonical filtration} of the $\lg$-module $L(\lambda)$.
\end{definition}

It follows
\[ L_l(\lambda)=\{ x(v) : x\in \U_l(\lg) \} \subseteq L(\lambda)\]
is the $\lpv$-module spanned by the highest weight vector $v$ and the elements in $\U_l(\lg)\subseteq \U(\lg)$.

Let $\ln(V_\bullet)_{op})$ be the opposite radical of the flag $V_\bullet$. It follows we get a direct sum decomposition
$\lg\cong \ln(V_\bullet)_{op}\oplus \lp$ of vector spaces. Let $z_1,\ldots ,z_d$ be a basis for $\ln(V_\bullet)_{op}$.
Let
\[ B_l(\lambda)=\{ z_1^{a_1}\cdots z_d^{a_1}(v): 0\leq \sum a_i \leq l\}.\]

Recall that by Theorem 4.1 we may do the following: There is a flag
\[ V_\bullet: 0 \neq V_1 \subsetneq V_2 \subsetneq \cdots \subsetneq V_k \subsetneq V \]
with $dim_K(V_i)=n_i$ and integers $l_1,\ldots ,l_k\geq 1$ with $\lambda=\sum_{i=1}^k l_i\omega_{n_i}$.
Theorem 4.1 gives an explicit construction of the highest weight vector $v$ and the irreducible finite dimensional
$\lg$-module $V(\lambda)$. Let $L_v$ be the line in $V(\lambda)$ spanned by $v$ and let
\[ \rho_v:\lp \rightarrow \End(L_v)\]
be the character associated to $L_v$. Let $char(\rho_v)$ be the character ideal of the character $\rho_v$.
We may choose a set of positive roots $R_+$ for the pair $(\lg, \lh)$ as follows:
\[ R_+=\{ L_i-L_j : 1\leq i < j \leq n\}.\]
It follows we get a Cartan decomposition $\lg=\ln_{-}\oplus \lh \oplus \ln_{+}$
where
\[ \ln_{+}=\oplus_{\alpha \in R_+}\lg_\alpha.\]
Let
\[ B=\{L_1-L_2,L_2-L_3,\ldots ,L_{n-1}-L_n\} \]
be a basis for the roots $R$ and let $\beta_i=L_i-L_{i+1}$. Let $\lb_{+}=\lg\oplus \ln_{+}$ and let
$\lb_{-}=\lh \oplus \ln_{-}$. Let furthermore $\delta=\frac{1}{2}\sum_{\alpha\in R_{+}}\alpha$. It follows by \cite{dixmier}, Section 7 that 
there is an isomorphism $V(\lambda)\cong L(\lambda+\delta)$ of $\lg$-modules.
It follows $\lg_{\beta_i}=KE_{i,i+1}$ and $\lg_{-\beta_i}=KE_{i+1,i}$.
We get
\[  [E_{i,i}, E_{j,j}]=E_{i,i}-E_{j,j}.\]
Let $X_{-\beta_i}=E_{i+1,i} \in \lg_{-\beta_i}$ and $H_{\beta_i}=E_{i,i}-E_{i+1,i+1}\in [\lg_{\beta_i},\lg_{-\beta_i}]$.
Let $m_{\beta_i}=\lambda(H_{\beta_i})+1$.

\begin{lemma}The following holds: $m_{\beta_i}=l_j+1$ if $i=n_j$. $m_{\beta_j}=1$ if $i\neq n_j$.
\end{lemma}
\begin{proof} The proof is left to the reader.
\end{proof}

Let $m(\lambda)=min_i\{ l_1,\ldots ,l_k\}$.

\begin{lemma}\label{mainlemma}  Let $1\leq l \leq m(\lambda)$.
There is an equality $ann_l(v)=char_l(\rho_v)$ of sub-vector spaces in $\U(\lg)$.
\end{lemma}
\begin{proof} Let
\[ I(v)=\U(\lg)\ln_{+} + \sum_{x\in \lh} \U(\lg)(x-\lambda(x)1_{\lg}) .\]
By \cite{dixmier}, Proposition 7.2.7 it follows
\[ ann(v)=I(v)+\sum_{\beta\in B}\U(\ln_{-})X^{\m_\beta}_{-\beta}.\]
Let 
\[ I_l(v)=I(v)\cap \U_l(\lg) \]
and
\[ J_l(v)=(\sum_{\beta\in B} \U(\ln_{-})X^{m_\beta}_{-\beta})\cap \U_l(\lg).\]
It follows
\[ J_l(v)=\sum_{i\neq n_j}\U_{l-1}(\ln_{-})X_{-\beta_i}+\sum_{i=n_j}\U_{l-l_j-1}(\ln_{-})X^{l_j+1}_{-\beta_{n_j}}.\]
When $1\leq l \leq m(\lambda)$ it follows
\[ J_l(v)=\sum_{i\neq n_j}\U_{l-1}(\ln_{-})X_{-\beta_i}.\]
By definition
\[ char_l(\rho_v)\subseteq ann_l(v)\]
for all $l\geq 1$.
By the above argument we see that
\[ ann_l(v)=I_l(v)+J_l(v)\subseteq char_l(\rho_v)\]
when $1\leq l \leq m(\lambda)$. The Lemma follows.
\end{proof}

In the following Theorem we give a basis for $L_l(\lambda)$ and calculate its dimension.

\begin{theorem} \label{maintheorem} Let $1\leq l \leq m(\lambda)$.
\begin{align}
&\label{one}\text{The set $B_l(\lambda)$ is a basis for $L_l(\lambda)$.}\\
&\label{two}dim_K(L_l(\lambda))=\binom{\sum_{1\leq j<i\leq k+1}l_il_j +l}{\sum_{1\leq j<i\leq k+1}l_il_j}.
\end{align}
\end{theorem}
\begin{proof} There is by Proposition \ref{directsum} a direct sum decomposition
\[ \U_l(\lg)=\U_l(\ln(V_\bullet)_{op})\oplus char_l(\rho) \]
for all integers $l\geq 1$. Here $\rho$ is the character $\rho:\lp\rightarrow \End(L_v)$.
There is by Lemma \ref{mainlemma}  an equality
\[ ann_l(v)=char_l(\rho) \]
for $1\leq l \leq m(\lambda)$.
We get isomorphisms
\[ L_l(\lambda)\cong M_l(\lambda)/ann_l(v)\otimes L_v\cong \U_l(\ln(V_\bullet)_{op})\otimes L_v  \]
hence claim \ref{one} follows from the PBW theorem.
By Proposition \ref{dimension} and the PBW theorem we get
\[ dim_K(L_l(\lambda))=dim_K(\U_l(\ln(V_\bullet)_{op})=\binom{\sum_{1\leq j<i\leq k+1}l_il_j+l}{\sum_{1\leq j<i\leq k+1}l_il_j},\]
and claim \ref{two} follows. The theorem is proved.
\end{proof}

We get for any finite dimensional irreducible $\sl(E)$-module $L(\lambda)$ a basis and a dimension formula for all the terms in the canonical 
filtration
\[ L_1(\lambda)\subseteq L_2(\lambda) \subseteq \cdots \subseteq L_l(\lambda) \subseteq L(\lambda) .\]

\begin{example} The case $\sl(K^2)$. \end{example}

Let $V=K^2$ and $\lg=\sl(V)$. It follows the $\lg$-modules $\Sym^l(V)$ for $l\geq1$ are all irreducible finite dimensional
$\lg$-modules. If $V$ has basis $e_1,e_2$ it follows $\Sym^l(V)$ has highest weight vector $v=e_1^l$. We get the following calculation:
\[ ann_k(v)=char_k(\rho) \]
for $1\leq k \leq l-1$.
If $k\geq l$ it follows
\[ ann_k(v)=K\{y^l,y^{l+1}\ldots , y^k\}\oplus char_k(\rho).\]
Here $y\in \lg$ is the following matrix
\[
y=
\begin{pmatrix} 0 & 0 \\
                1 & 0 
\end{pmatrix}
.\]
There is an isomorphism of vector spaces
\[ \U(\lg)\cong K[y]\oplus char(\rho) \]
where $K[y]\cong \U(\ln_-)$ and $\ln_-$ is the abelian sub algebra of $\lg$ generated by $y$.

\begin{example} Symmetric powers of the standard $\sl(V)$-module.\end{example}

Let $V=K\{e_1,\ldots ,e_n\}$ and let $V_l=\Sym^l(V)$. It follows the vector $v=e_1^l$ is a highest weight vector for the irreducible
$\sl(V)$-module $V_l$ for all $l\geq 1$. It has highest weight $\lambda=lL_1$. Let $L(\lambda)=V_l$.
It follows from Theorem \ref{maintheorem} that 
\[ dim_K(L_k(\lambda))=\binom{n-1+k}{n-1} \]
for $1\leq k \leq l$.
We moreover get
\[ dim_K(L_l(\lambda))=\binom{n-1+l}{n-1}=dim_K(L(\lambda)).\]
We get an equality $L_l(\lambda)=L(\lambda)$ of $\lp$-modules.

\begin{example} The adjoint representation. \end{example}

Let 
\[ an: \sl(V) \rightarrow \End(\sl(V)) \]
be the adjoint representation. This representation is irreducible with highest weight $\lambda=\omega_1+\omega_{n-1}$. 
One checks we get a strict inclusion
\[ L_1(\lambda)\subsetneq L(\lambda) \]
and an equality
\[ L_2(\lambda)=L(\lambda) \]
of $\sl(V)$-modules.

\begin{example} A geometric construction of the canonical filtration $L_l(\lambda)$ .\end{example}

The canonical filtration $L_l(\lambda)$ is a filtration of $\lp(V_\bullet)$-modules where $\lp(V_\bullet)$ in $\sl(V)$ is the stabilizer Lie algebra
of a flag $V_\bullet$ in $V$. There is a parabolic subgroup $P$ in $\SL(V)$ with $Lie(P)=\lp(V_\bullet)$ and the terms $L_l(\lambda)$ are $P$-modules.
There is an equivalence of categories between the category of finite dimensional $P$-modules and the category of locally free finite rank
$\O_{\SL(V)/P}$-modules with an $\SL(V)$-linearization. Hence to any term $L_l(\lambda)$ in the canonical filtration one gets a vector bundle
$J_l(\lambda)$ on the flag variety $\SL(V)/P$. In the preprint \cite{maakestad0} we give for any irreducible finite dimensional $\SL(V)$-module 
$L(\lambda)$ a geometric construction of the vector bundle $J_l(\lambda)$
in terms of jet bundles of line bundles on $\SL(V)/P$. The proof use higher direct images of sheaves on $\SL(V)/P$, Taylor maps and general properties
of jet bundles. It also use the dimension formula \ref{two} from Theorem \ref{maintheorem}.

\begin{example} \label{discriminants} Applications: Discriminants of linear systems on flag varieties.\end{example}

Let $G$ be a linearly reductive  algebraic group over an algebraically closed field $K$ of characteristic zero and let
$P$ be a parabolic subgroup of $G$. Let $\L$ be a $G$-linearized line bundle on $G/P$ and let $J^l(\L)$ be the $l$'th
jet bundle of $\L$. There is by \cite{maakestad2}, Section 4 a resolution
\[ 0 \rightarrow \O(-r)_Y\otimes \wedge^r J^l(\L)^*_Y \rightarrow \cdots \rightarrow \O(-2)_Y\otimes \wedge^2J^l(\L)^*_Y \rightarrow \]
\[ \O(-1)_Y\otimes J^l(\L)^*_Y \rightarrow \O_Y\rightarrow \O_{I_l(\L)}\rightarrow 0 \]
where $I_l(\L)$ is the $l$'th incidence scheme of the line bundle $\L$. There is a canonical map
\[ p: I_l(\L)\rightarrow D_l(\L) \]
where $D_l(\L)$ is the $l$'th discriminant of $\L$. The higher direct image of the sheaf $\O(-i)_Y\otimes \wedge^iJ^l(\L)^*_Y$ equals
\begin{align}
&\label{coho}  \O(-i)\otimes \H^j(G/P, \wedge^i J_l(\L)^*).
\end{align}
Let $L(\lambda)=H^0(G/P, \L)^*$.  It follows by the results of \cite{maakestad0} that there is an isomorphism
\[ L_l(\lambda)\cong J^l(\L)(e)^*\]
of $P$-modules. Hence the cohomology group $\H^j(G/P, \wedge^i J_l(\L)^*)$ may be calculated using Bott's Theorem on the calculation 
of the cohomology of bundles on the flag variety. This might give information on the syzygies of the discriminant $D_l(\L)$ as it does in the case
of the syzygies of determinantal varieties as done by Lascoux.

\begin{example} Canonical bases for semi simple algebraic groups.\end{example}

Let $G$ be a semi simple linear algebraic group and $L$ a finite dimensional irreducible $G$-module with highest weight $\lambda$
and highest weight vector $v\in L$. Much work has been devoted to the construction of canonical bases in $L$. In the paper \cite{lusztig} 
the author constructs a \emph{universal} canonical basis $B$ in a quantized enveloping algebra $U$ associated to a root system. He shows that 
the basis $B$ specialize to a basis for $L$ for all highest weights $\lambda$. One could speculate about the existence of a universal canonical
filtration with a universal canonical basis $B_l$ specializing to a basis for $L_l$ for every highest weight $\lambda$. One should then compare
this basis to the basis given in Theorem \ref{maintheorem} for every highest weight $\lambda$. 

\begin{example} Standard monomial theory.\end{example}

Let $G$ be a semi simple linear algebraic group over an algebraically closed field of characteristic zero.
In the paper \cite{brion} the authors construct for every finite dimensional irreducible $G$-module $L$ a basis using a geometric 
approach. By the Borel-Weil-Bott formula one may realize the $G$-module $L$ as the global sections of a line bundle $\L$ on a quotient
$G/P$ where $P$ is a parabolic subgroup of $G$. They embed $G/P$ into a product of projective spaces
\[ i: G/P \rightarrow \P^{N_1}\times \cdots \times \P^{N_k} \]
and use \emph{standard monomials} as a basis for $L$. These monomials are defined in terms of the embedding $i$.
It is not clear how the basis constructed in \cite{brion} is related to the basis constructed in Theorem \ref{maintheorem}.

\section{Appendix: Irreducible finite dimensional $\sl(V)$-modules}

Let in this section $\omega_i=L_1+\cdots +L_i$ for $i=1,\ldots , n-1$. It follows $\omega_1,\ldots ,\omega_{N-1}$ are the 
\emph{fundamental weights} for the Lie algebra $\lg=\sl(V)$. Let $\lg=\ln_{-}\oplus \lh \oplus \ln_{+}$ be the standard Cartan decomposition where
$\lh$ is the abelian sub algebra of diagonal matrices with trace zero.
In Dixmier's book \cite{dixmier} the following classification result 
is proved: Let $V$ be a finite dimensional irreducible $\lg$-module. It follows there is a unique line $L_v$ in $V$ with generator $v$ 
satisfying the following properties: $V$ is generated by the elements $x(v)$ for $x\in \ln_{-}$. For all $x\in \lh$ it follows 
$x(v)=\lambda(x)v$ where $\lambda=\sum_{i=1}^{N-1}\alpha_i \omega_i$ and where $\alpha_i\geq 0$ are integers. Moreover $x(v)=0$ for all 
$x\in \ln_{+}$. We say $\lambda\in \lh^*$ is the \emph{highest weight} for the module $V$. The vector $v$ is called a \emph{highest weight vector} 
for $V$. We write $L(\lambda)$ for the module $V$. 

By the general theory of \cite{dixmier} it follows $L(\lambda)$ is classified by
the highest weight $\lambda$:

Let $1\leq d_1 < d_2 < \cdots < d_k \leq N-1$ and let $l_1,\ldots l_k\geq 1$ be integers. Let $\lambda=\sum_{i=1}^k l_i\omega_{d_i}\in \lh^*$.
There is a unique irreducible finite dimensional $\sl(V)$-module $L(\lambda)$ with highest weight $\lambda$.
Moreover every irreducible finite dimensional $\sl(V)$-module arise in this way. Let $\underline{l}=(l_1,\ldots ,l_k)$ and
$\underline{n}=(n_1,\ldots ,n_k)$. Let $E=K\{e_1,\ldots ,e_N\}$ and let $V_i=K\{e_1,\ldots ,e_{n_i}\}$. It follows we get a flag
\[ V_\bullet: 0\neq V_1\subsetneq V_2\subsetneq \cdots \subsetneq V_{k+1}=V \]
of vector spaces in $V$. Let
\[ W(\underline{l},\underline{n})=\Sym^{l_1}(\wedge^{n_1}E)\otimes \cdots \Sym^{l_k}(\wedge^{n_k}E) .\]
It follows $W(\underline{l},\underline{n})$ is an $\sl(E)$-module. Let  $w_i=\wedge^{n_i}V_i$ and let 
\[ v=\Sym^{l_1}(w_1)\otimes \cdots \otimes \Sym^{l_k}(w_k) \in W(\underline{l},\underline{n}).\]
One checks the following: For all $x\in \ln_{+}$ it follows $x(v)=0$. For all $x\in \lh$ it follows $x(v)=\lambda(x)v$.
Let $\U(\lg)v$ be the sub $\lg$-module generated by $v$.

The following Theorem gives an explicit construction of all irreducible finite dimensional $\lg$-modules $L(\lambda)$. 
It also constructs the stabilizer Lie algebra of a highest weight vector $v$ in $L(\lambda)$ in terms of the stabilizer Lie algebra
of a flag $V_\bullet$.

\begin{theorem} \label{irreducible}  The following holds: $\U(\lg)v\cong L(\lambda)$. There is an equality
\[ \lp(V_\bullet)=\{ x\in \lg : x(v)\in L_v\}.\]
Hence the stabilizer Lie algebra of a highest weight vector $v$ in $L(\lambda)$ equals the stabilizer Lie algebra of the flag $V_\bullet$.
\end{theorem}
\begin{proof} The proof is left to the reader as an exercise.
\end{proof}

Note: The classification of irreducible finite dimensional $\sl(V)$-modules is the same as the classification of finite dimensional
irreducible $\SL(V)$-modules, where $\SL(V)$ is the linear algebraic group of linear automorphisms of $V$ with determinant one.
One may formulate the notion of a highest weight vector and highest weight in terms of algebraic groups and one gets similar theorems
for algebraic groups as in the theory of Lie algebras.

\textbf{Acknowledgements} Thanks to Michel Brion, Alexei Roudakov, David Vogan and an anonymous referee for comments on 
the contents of this paper.

\end{document}